# SAMPLING THE LOWEST EIGENFUNCTION TO RECOVER THE POTENTIAL IN A ONE–DIMENSIONAL SCHRÖDINGER EQUATION

ROB RAHM †

Abstract. We consider the BVP $-y'' + qy = \lambda y$ with $y(0) = y(1) = 0$. The inverse spectral problems asks one to recover $q$ from spectral information. In this paper, we present a very simple method to recover a potential by sampling one eigenfunction. The spectral asymptotics imply that for larger modes, more and more information is lost due to imprecise measurements (i.e. relative errors *increases*) and so it is advantageous to use data from lower modes. Our method also allows us to recover "any" potential from *one* boundary condition.

## 1. Introduction

The goal of this paper is to present a method for recovering the potential in a one–dimensional Schrödinger equation by sampling only the lowest mode eigenfunction. In particular, if $\{\lambda_k\}_{k=1}^{\infty}$ and $\{y(x, \lambda_k; q)\}_{k=1}^{\infty}$ are the eigenvalues (in increasing order) and eigenfunctions for

$$-y'' + qy = \lambda y, \qquad\qquad y(0) = y(1) = 0, \tag{1.1}$$

then we want to recover $q$ by sampling $y(x, \lambda_1; q)$ at a relatively small number of points.

This equation (the Schrödinger equation) is ubiquitous both as a model problem and as a sort of canonical form of any second order linear ODE. Typically, one wants to determine the potential $q$ from some (finite) spectral data. For example, in the context of strings, $q$ is related to the density of the string via the Liouville transform and one would like to recover the density by measuring properties of the resonances and corresponding eigenmodes.

The equation $-y'' + qy = \lambda y$ is a perturbation of the equation $-y'' = \lambda y$ and as $k \to \infty$ the eigenfunctions and eigenvalues all converge (somewhat rapidly) to the spectral data in the $q \equiv 0$ case. Thus, in a sense, there is "more" usable information in the lower modes than in the higher modes. So it is desirable to use only the lower modes to recover $q$ (this is discussed more below).

Our method is based on derivative formulas from [8]. In particular, for a fixed $\lambda$ and $x$ we consider the map:

$$q \mapsto y_2(x, \lambda, q).$$

(Throughout, $y_1(x, \lambda, q)$ and $y_2(x, \lambda, q)$ will denote the solutions to $-y'' + qy = \lambda q$ with initial conditions $y_1(0, \lambda, q) = 1$ and $y_1'(0, \lambda, q) = 0$ and $y_2(0, \lambda, q) = 0$ and $y_2'(0, \lambda, q) = 1$.) The derivative of this function is a linear map from $L^2$ to $\mathbb{C}$. By the Riesz Representation Theorem, this can be given by an integral. In [8], it was shown that the derivative is:

$$\frac{\partial y_2(x, \lambda, q)}{\partial q}(v) = \int_{t=0}^{1} K(t, x, \lambda, q) v(t) dt,$$

† robrahm@tamu.edu.

where:

$$K(t,x,\lambda,q) = y_2(t,\lambda,q)\left(y_1(t,\lambda,q)y_2(x,\lambda,q) - y_1(x,\lambda,q)y_2(t,\lambda,q)\right)\mathbb{1}_{[0,x]}(t).$$

To apply this, we will consider a restricted map from a finite dimension subspace of $L^2$ (i.e. we will consider $q$ to be a finite linear combination of basis functions). Each sampling of an eigenfunction corresponds to one equation. So, if we have $n$ sample points, we should be able to solve the equation for $q$ in a $n$–dimensional subspace.

The paper is organized in the following way. In Section 2 we discuss some previous work in this area (there has been a lot). In Section 3 we give our inverse method and illustrate it with some examples. Finally, in Section 4 we give the code listings used here.

## 2. Previous Work and Current Work

The inverse problem that we discuss here is old and well-researched. We can't discuss all previous work, but we discuss some of the results that we find to be most relevant.

In 1946, Borg proves in [1] that if two spectra are given (corresponding to two different boundary conditions) then, in principle, one can recover $q$. Levinson proved in [5] (with a different and much shorter proof) that one can also recover the boundary condition. In 1951, Gel'fand and Levitan ([2]) show that the potential can be recovered from a spectrum and some information on the eigenfunctions (in particular the norming constants).

In practice, to recover the potential, we are only given a finite set of spectral data. So precisely speaking there isn't a uniqueness theorem. So any algorithm will necessarily only produce a potential that approximates the target potential. In [3], Hald gives a method to recover an even potential from a spectrum and boundary conditions. A result similar to ours by Hald and McLaughlin in [4] gives a way to determine a potential from nodal data of the eigenfunctions. Another method that is similar to ours is in [6] and, for example, [10].

As discussed above, the spectral data for the Schrödinger operator satisfy the following asymptotics: (see [8]):

$$\lambda_k(q) = (k\pi)^2 + \int q(x)\,dx + \alpha_k(q), \qquad y_2(x,\lambda;q) = \frac{\sin\sqrt{\lambda}x}{\sqrt{\lambda}} + O(\frac{1}{\lambda}),$$

where $\{\alpha_k\}$ is an $\ell^2$ sequence. In other words, spectral data is "centered" around (and converges to) the corresponding data for the constant potential equal to the average value of $q$. For example, if we consider just well–behaved, mean–value zero potentials, then

$$\lambda_{10}(q) = (10\pi)^2 + \alpha_k(q) \approx 987 + \alpha_k(q).$$

Typically, we can take $\alpha_k$ to be bounded in absolute value by $k^{-1}$. Therefore, in the tenth eigenvalue, the information that distinguishes one potential from another is contained in an interval of radius $\frac{1}{10}$ centered at 987. In other words, for the tenth eigenvalue to be usable, a relative precision of (at least) .01 *percent* is needed.

So, it is better to use lower–mode data. Such techniques were also considered in, for example, [7] where the lowest eigenvalue corresponding to different boundary conditions is used and in [9] where (in essence) known masses are added to the potential and the lowest eigenvalue is measured. The current paper uses sampling data from the lowest mode eigenfunction – though as is shown below, we can also use higher modes.

Finally, we mention that in addition to the sample data, we need to know both the eigenvalue and the value of $y'(0)$. So in a sense, we can approximate the projection of $q$ onto the first $n$ basis functions by collecting $n+2$ pieces of data.

**Acknowledgment.** I'd like to thank Bill Rundell for helpful discussions and information regarding the topic of this paper.

## 3. The Method

In this section, we describe our method and provide a discussion of its efficacy - mentioning both its advantages and disadvantages. Psychologically, We concentrate our attention to the setting of Dirichlet–Dirichlet boundary conditions since this corresponds to (for example) the setting of a clamped string. However, the methods we use only really require a Dirichlet condition at the left endpoint. After a discussion of the method, we present some graphs of our algorithm's reconstruction of some potentials.

3.1. **Description of the Method.** Fix $\lambda$ and define an operator $T : L^2([0,1]) \to C^2([0,1])$ by:

$$Tq(x) = y_2(x, \lambda, q)$$

where $y_2(x, \lambda, q)$ is the solution to $-y'' + py = \lambda y$ with $y(0) = 0$ and $y'(0) = 1$. Let $y_1$ be the solution with $y(0) = 1$ and $y'(0) = 0$ (this is the notation used in [8]). It was shown in [8] that the derivative of this operator (with respect to the potential) is:

$$T'_p(v) = \int_{t=0}^{x} y_2(t)\left(y_1(t)y_2(x) - y_1(x)y_2(t)\right)v(t)\,dt =: \int_{t=0}^{x} K_p(t, x; \lambda)v(t)\,dt.$$

If $q$ is near $p$, linearization then gives:

$$y_2(x, \lambda, q) \simeq y_2(x, \lambda, p) + \int_{t=0}^{x} K_p(t, x; \lambda)(q(t) - p(t))\,dt. \tag{3.1}$$

In practice, we will use a quasi–Newton method. In this case, the kernel $K_p$ is replaced by the kernel $K_0$ which is $K(x, t) = s_\lambda(t)s_\lambda(x - t)$ where $s_\lambda(x) := \frac{\sin\sqrt{\lambda}x}{\sqrt{\lambda}}$. This gives the quasi–Newton iteration:

$$q_0 = 0, \quad q_{k+1} = q_k + \delta, \quad \text{where} \quad y_2(x, \lambda, q) = y(x, \lambda, q_k) + \int_{t=0}^{x} s_\lambda(t)s_\lambda(x - t)(t, x; \lambda)\delta(t)\,dt.$$

To solve this equation for $\delta$, we assume that $q$ is a linear combination of some basis functions $\{\varphi_l\}$ (e.g. $\varphi_l(x) = \cos(2(l-1)\pi x)$ if we assume *a priori* that $q$ is even.) Then $\delta$ is a sum of those same basis functions: $\delta(x) = \sum_{l=1}^{n} d_l \varphi_l(x)$. Since there are $n$ unknowns, we need to convert our spectral data into $n$ equations.

For example, assume that we know $y_2(x, \lambda, q)$ evaluated at $n$ points - say $x_1, \ldots, x_n$. The quasi–Newton method gives the system:

$$\begin{pmatrix} y_2(x_1, \lambda, q) \\ \vdots \\ y_2(x_n, \lambda, q) \end{pmatrix} = \begin{pmatrix} y_2(x_1, \lambda, q_k) \\ \vdots \\ y_2(x_1, \lambda, q_k) \end{pmatrix} + \begin{pmatrix} \int_{t=0}^{x_1} s_\lambda(t)s_\lambda(x_1 - t)\delta(t)\,dt \\ \vdots \\ \int_{t=0}^{x_n} s_\lambda(t)s_\lambda(x_n - t)\delta(t)\,dt \end{pmatrix}. \tag{3.2}$$

This can be written as a matrix–vector equation. If $J$ is the matrix whose $j, l$ entry is:

$$\int_{t=0}^{x_j} s_\lambda(t)s_\lambda(1 - t)\varphi_l(t)\,dt,$$

and $\vec{y_2}(1,\lambda,q_k)$ is the first term, we can write this as:

$$\vec{0} = \vec{y_2}(1,\lambda,q_k) + J\vec{d},$$

where $\delta(t) = \sum_{l=1}^{n} d_l \varphi_l(t)$.

3.2. **Discussion of Recovery.** In this subsection, we explain why our algorithm is able to recover the potential.

Using the derivative formula and the mean value theorem, we can estimate:

$$|y_2(x,\lambda,q) - y_2(x,\lambda,p)| \leq \int_0^x |K_p(x,t)(q(t)-p(t)|\,dt + O(\|q-p\|^2).$$

In many cases, the kernel is bounded by 1. If $q(x) = \sum_{l=0}^{\infty} a_l \cos(l\pi k)$ and $q_N(x) = \sum_{l=0}^{N-1} a_l \cos(l\pi x)$ then standard estimates show that $\|q - q_N\| \leq \frac{1}{N}$.

In addition, by the injectivity of the Jacobian in the quasi–Newton iteration scheme, if $q_1$ and $q_2$ are linear combinations of the first N basis functions, and if their values at N sample points are close, then their difference $q_1 - q_2$ must be small (how small depends on precise quantitative properties of the Jacobian - and this depends on the choice of sample points which is discussed below.)

Putting these two facts together shows that the potential the algorithm recovers is close to projection of the target potential on the first N basis functions.

Recall that the Jacobian is the matrix whose $j, l$ entry is:

$$\int_{t=0}^{x_j} s_\lambda(t)s_\lambda(x_j - t)\varphi_l(t)dt, \qquad \text{where} \qquad s_\lambda(t) = \frac{\sin\sqrt{\lambda}t}{\sqrt{\lambda}}.$$

It is of course relevant to ask if this matrix is well–behaved. For example, is it injective and what is its condition number? These data depend on $\lambda$ and the sample points. Standard computations (i.e. explicitly computing singular values) shows that for evenly spaces sample points at $\sqrt{\lambda} = k\pi$ then this matrix is invective. Since $\lambda_k$ is close to $(k\pi)^2$ this indicates (via continuity of the determinant function), that for most potentials the Jacobian is injective. We have never found a potential where the Jacobian isn't invective and, in practice we use the Moore–Penrose pseudoinverse so that we can delete smaller singular values.

In practice, if we let $C(x_1, \ldots, x_n)$ be the condition number of this Jacobian with $\sqrt{\lambda} = k\pi$, then we can choose the sample points by minimizing this condition number. This is illustrated in the examples below.

A limitation of this method is that we must know a good basis with which to approximate q *a priori*. In our method, once we have sampled the eigenfunction, we can then run the algorithm using different basis functions for the $\varphi_l(x)$. With out knowing q, how can we determine which basis gives us a better approximation?

A related question is this: suppose we have two orthonormal bases $\{\varphi_l\}$ and $\{\psi_l\}$ and we know the inner products $\{\langle q, \varphi_l\rangle\}_{l=1}^{N}$ and $\{\langle q, \psi_l\rangle\}_{l=1}^{N}$ (i.e. we know the projection of q onto the first N basis functions.) We can determine which projection is a better $L^2$ approximation in the following way. Since

$$\|q\|_{L^2} = \sum_{l=1}^{\infty} (\langle q, \varphi_l\rangle)^2 = \sum_{l=1}^{\infty} (\langle q, \psi_l\rangle)^2$$

the better approximation is going to be the one that captures more of the norm, that is whichever projection has the larger $L^2$ norm. Since the algorithm produces an approximation to the projection of $q$ onto the basis $\{\varphi_l\}$, if we have two (or more) approximations to $q$, we just select the one that has the higher $L^2$ norm. This selection criteria is illustrated in the final example.

3.3. **Examples Using Sampling.** In this section, we give examples of our method. In addition to sampling the lowest mode eigenfunction, we also sample higher modes to illustrate the method further. We also use different bases.

For each potential and each sampled eigenfunction, we will plot the potential and its reconstruction on the plot on the left, and the plot on the right is the eigenfunction we sampled. The asterisks in the plot of the eigenfunction indicate the location of the sample points.

*Example* 1. The first example deals with $q(x) = 1 - \exp(-20(x - \frac{1}{2})^2)$. We also use the even cosine basis $\cos 2(k-1)\pi x$ for $k = 1, 2, \ldots$.

For the first reconstruction, we use the first eigenfunction and three sample points chosen by the optimization procedure explained above.

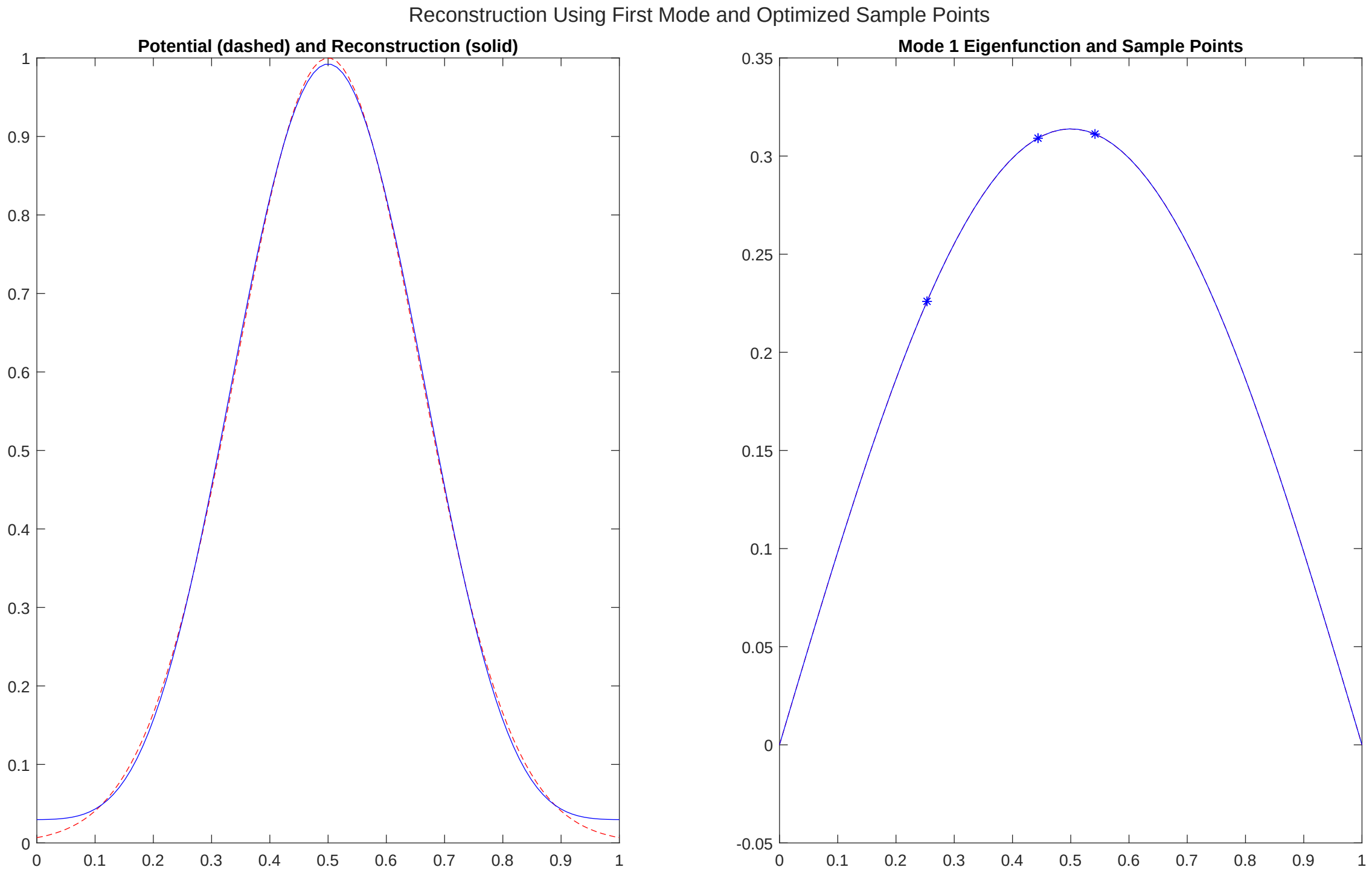

For the next reconstruction, we use the second eigenfunction and three sample points chosen by the optimization procedure explained above.

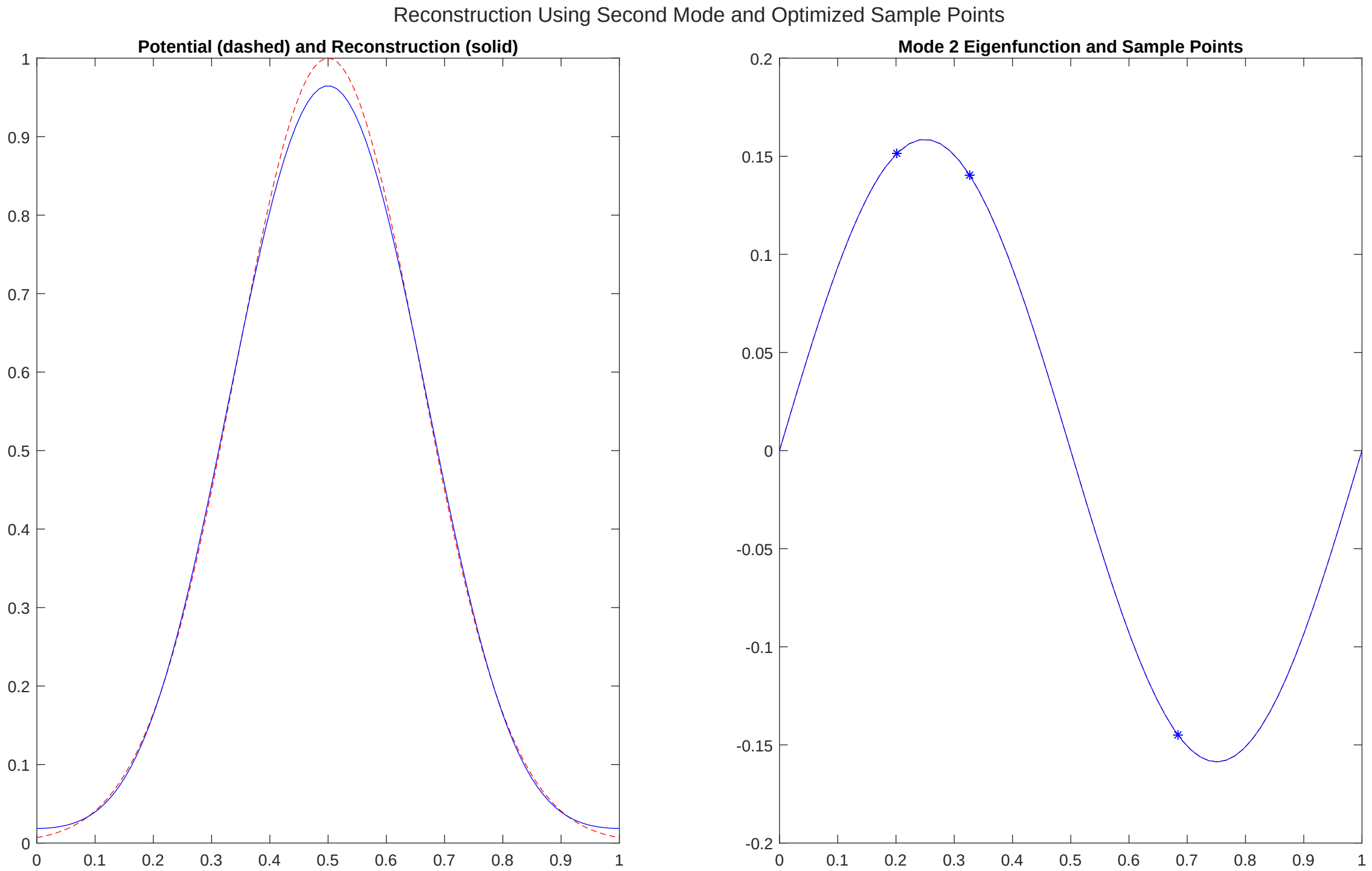

For the next reconstruction, we use the first eigenfunction and three equally–spaced sample points.

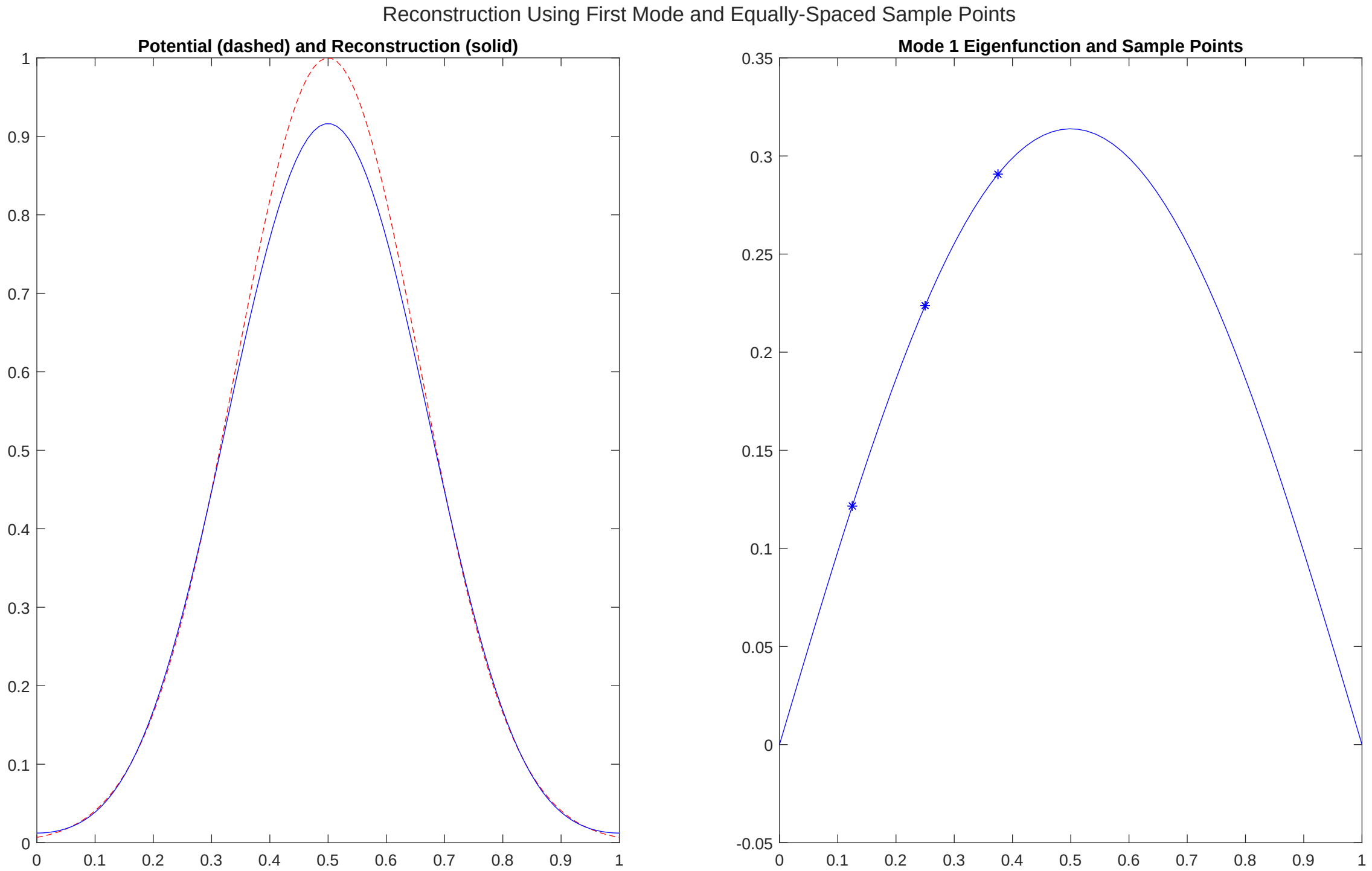

*Example* 2. For the next example, we use the function $q(x) = 1 - |x - .25|\, \mathbb{1}_{[0,.5]}(x) - |x - .75|\, \mathbb{1}_{[.5,1]}(x)$. We also use the even cosine basis $\cos 2(k-1)\pi x$ for $k = 1, 2, \ldots$.

In the first reconstruction, we use three optimally spaced sample points and the first eigenfunction.

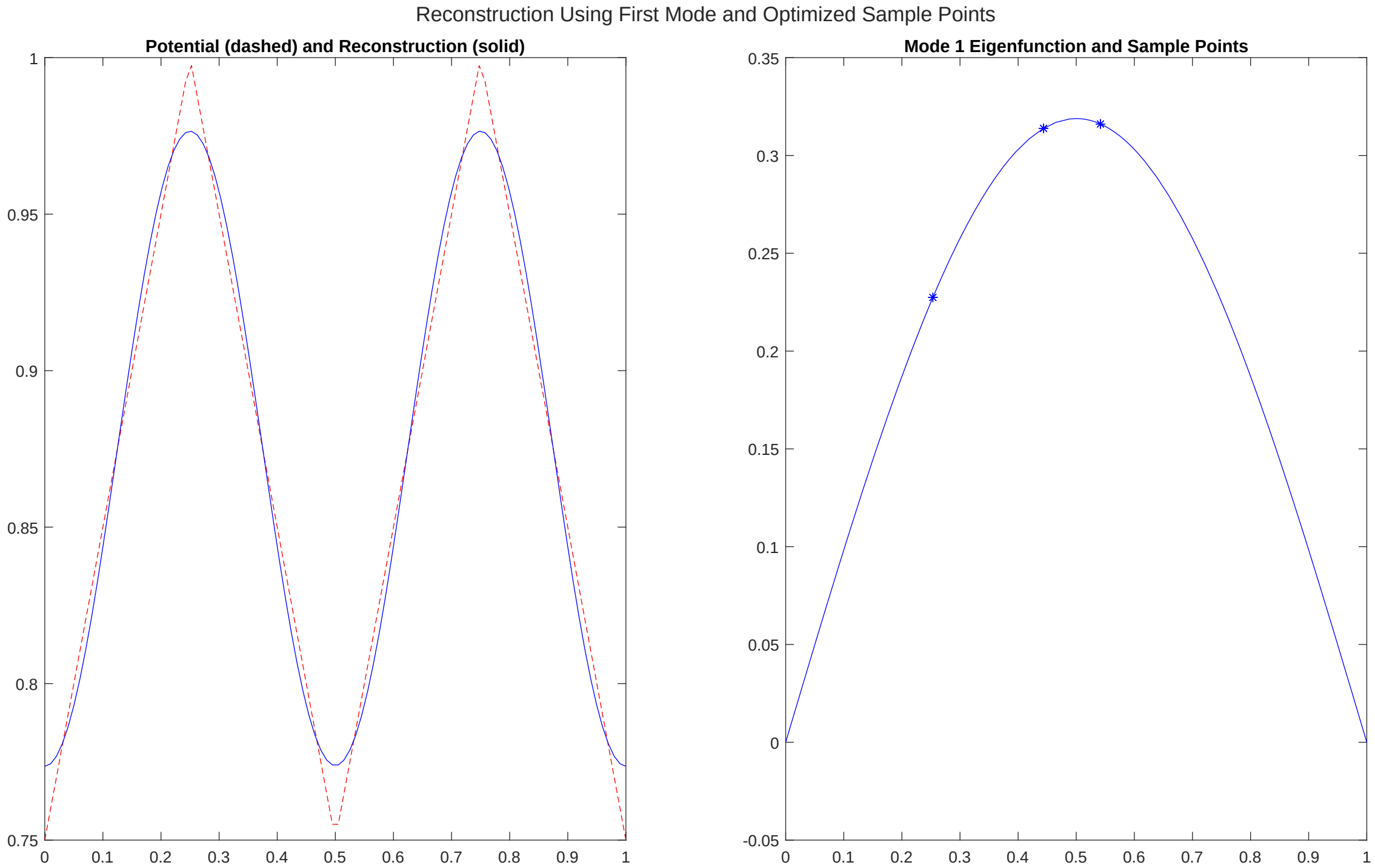

In the next reconstruction, we use three equally spaced sample points and the first eigenfunction.

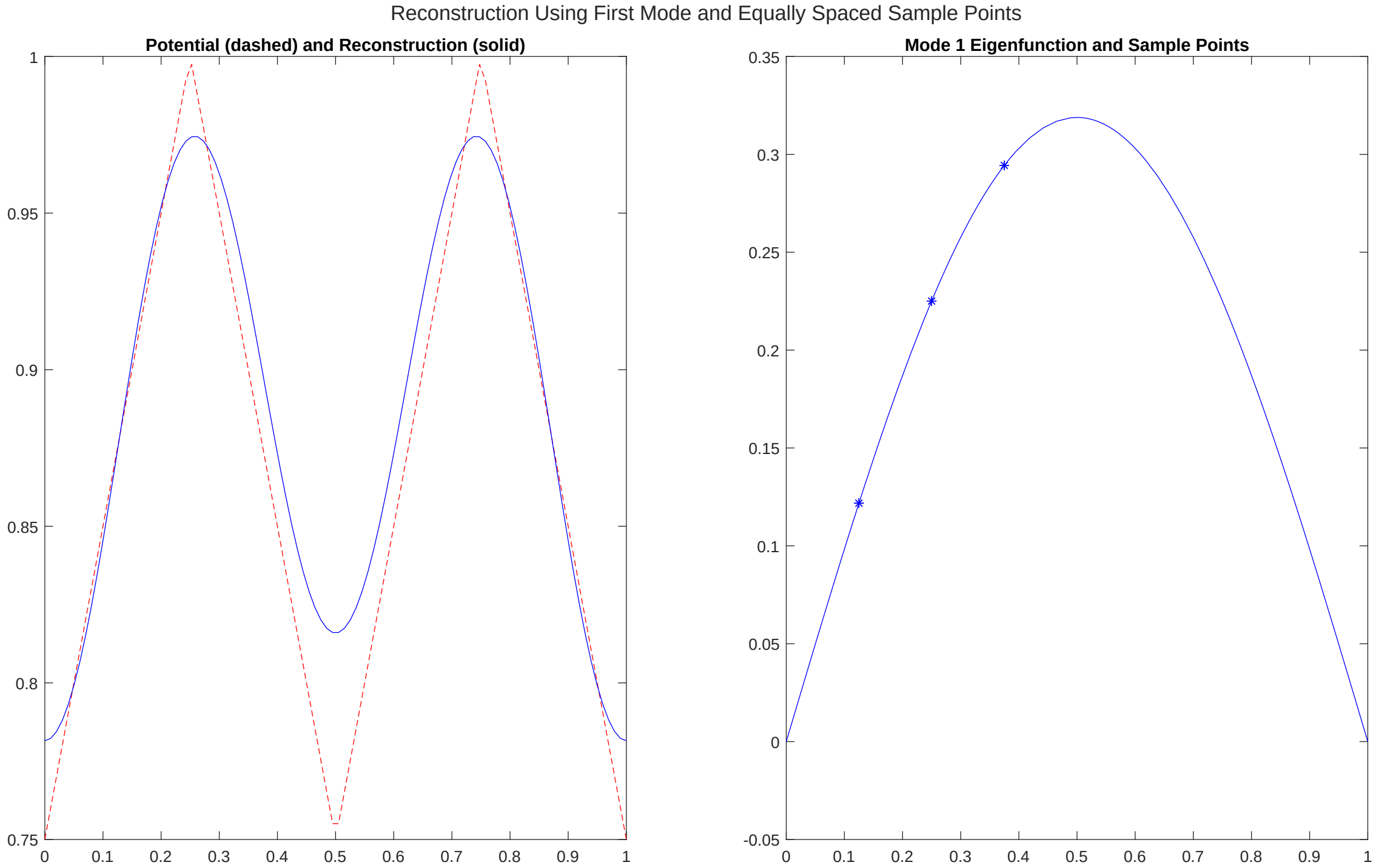

*Example* 3.

*Example* 4. In the next example, we work with the potential:

$$q(x) = p(x - .5), \qquad \text{where} \qquad p(x) = (2x)^6 - 3(2x)^4 + (2x)^2 - 1.$$

We use the even Legendre polynomial basis. In the first case, we use three optimally spaced sample points and the first eigenfunction:

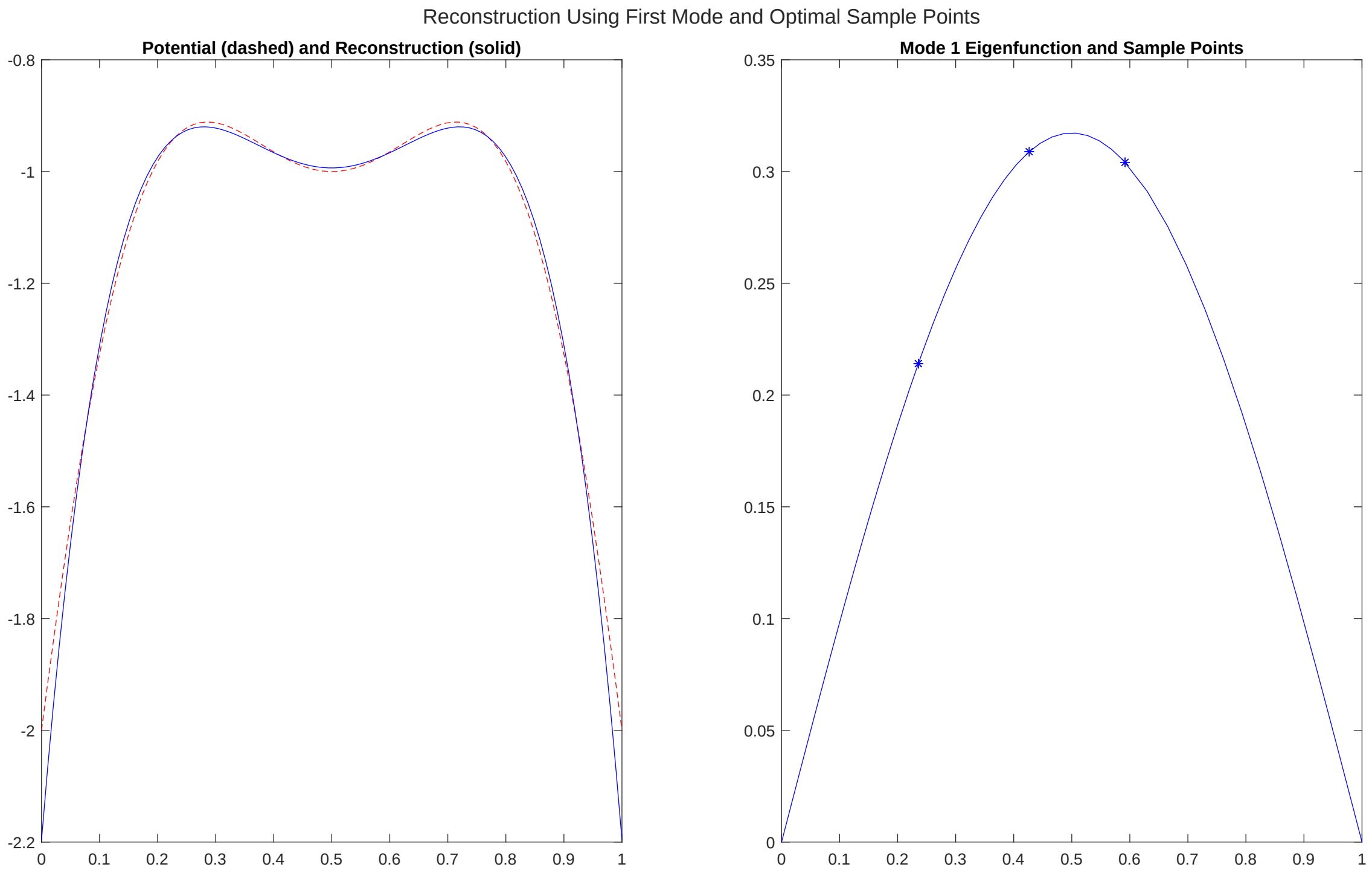

In the second case we use four optimally spaced points and the first eigenfunction. Of course, this means we are using degree four Legendre polynomials and so the recovery is perfect.

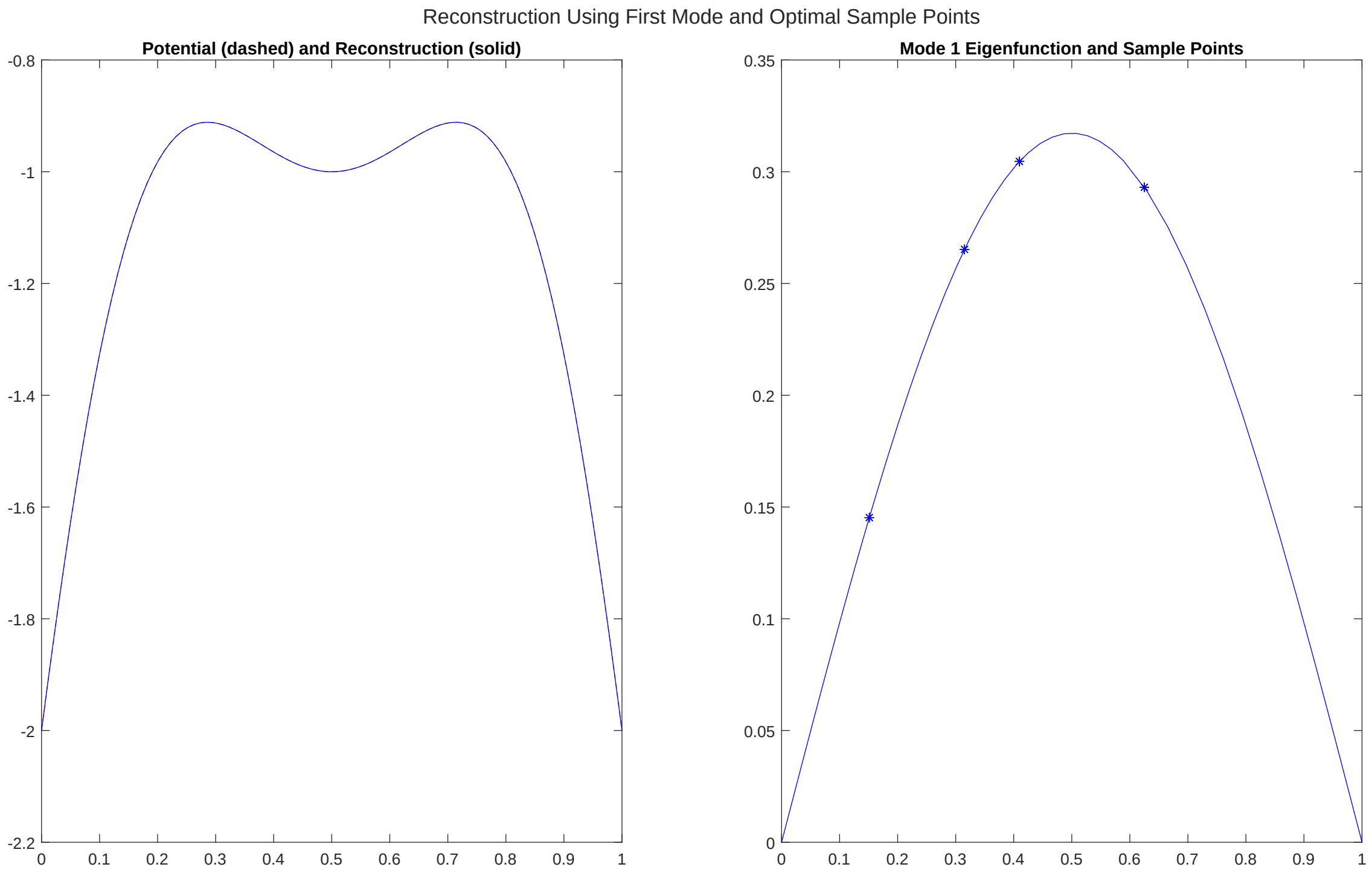

In the examples below we use equally spaced sample points:

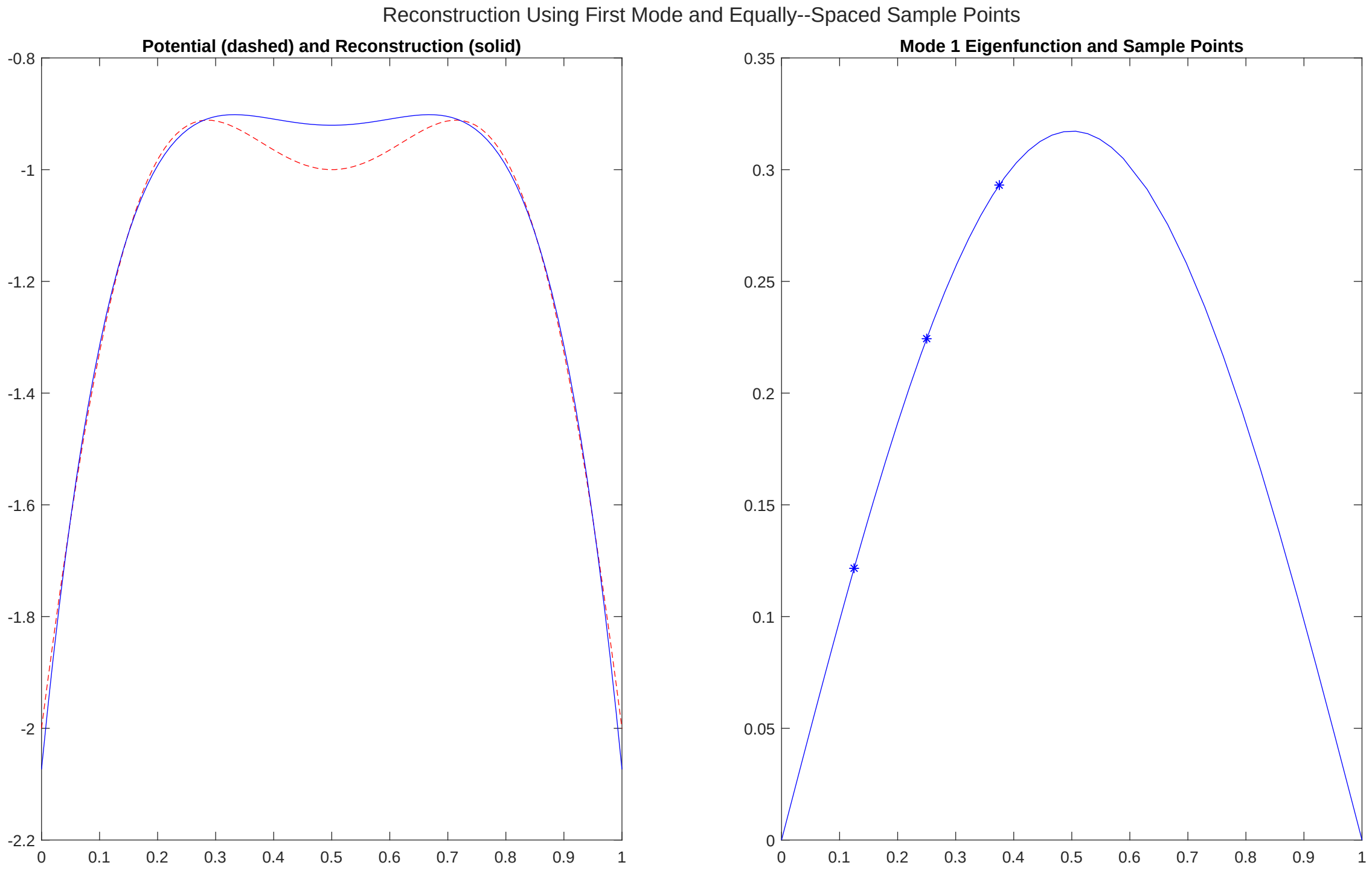

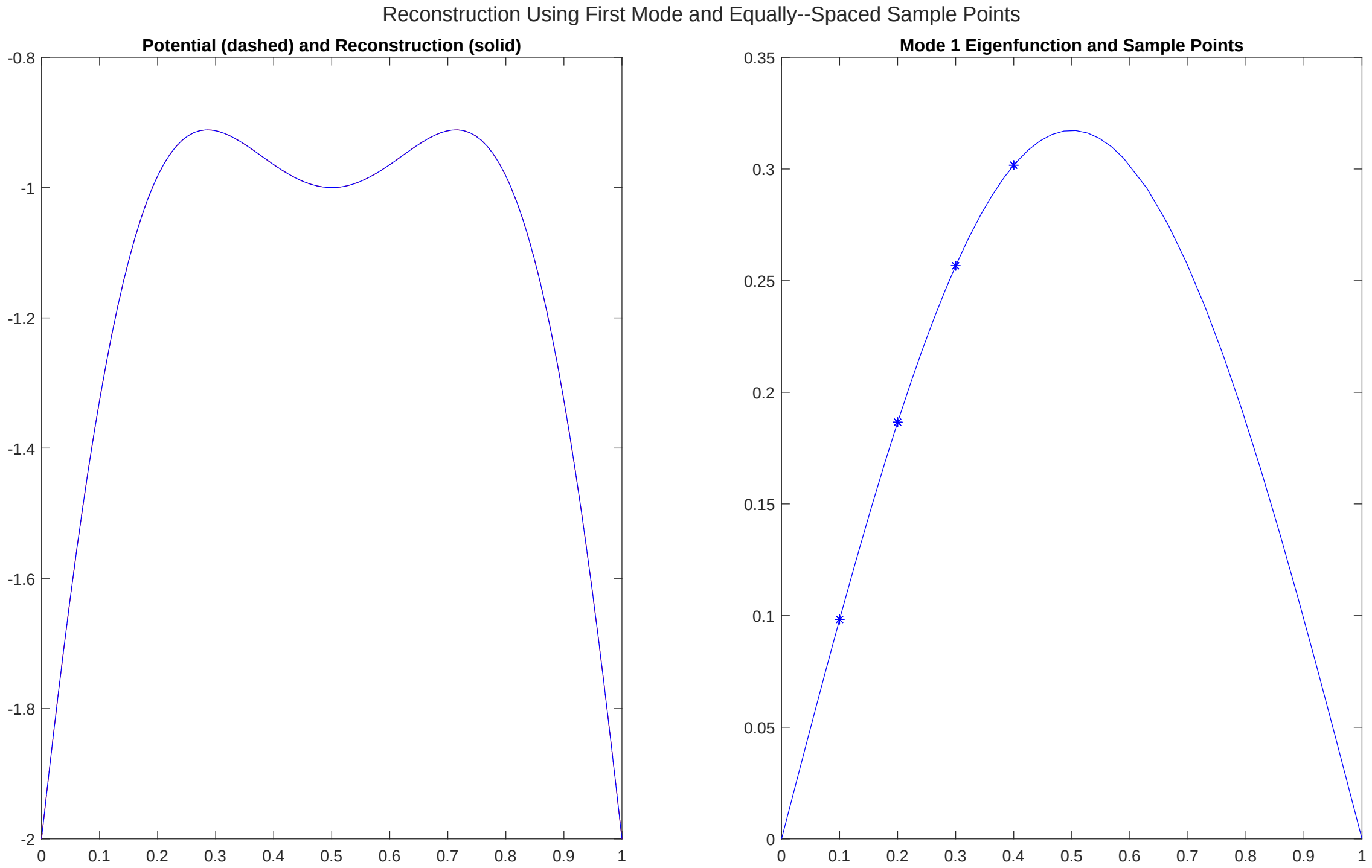

*Example* 5. Next we test our methods on the potential

$$q(t) = 1 + t + .3\cos(2\pi t) - .1\sin(2\pi t) + \cos(4\pi t) + .56\sin(4\pi t),$$

and we use the Legendre basis again with six sample points.

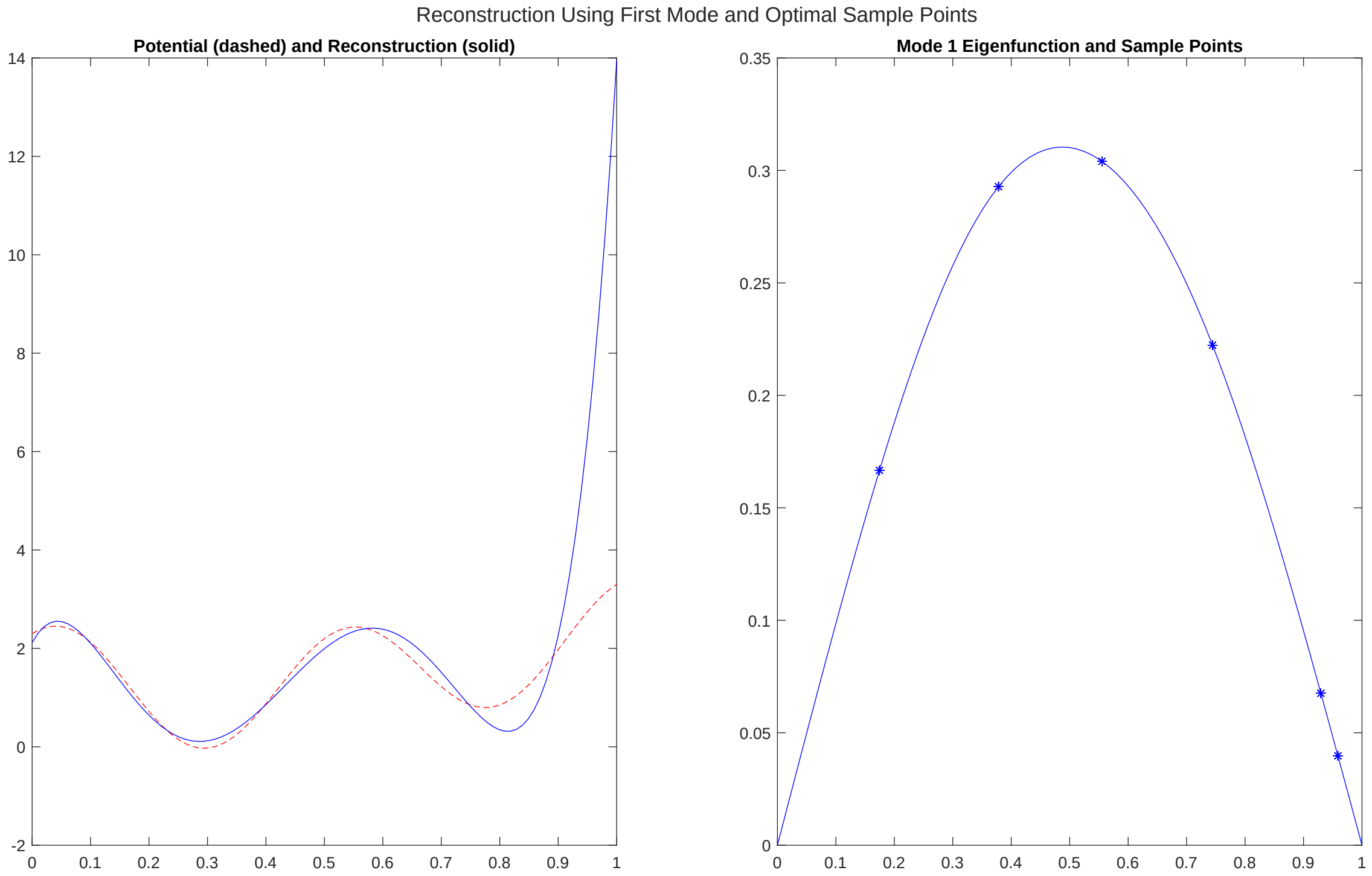

*Example* 6. Finally, we use the potential $q(t) = t^6 + t^5 - t$. Once again we use the Legendre basis. Unsurprisingly, the reconstruction is perfect.

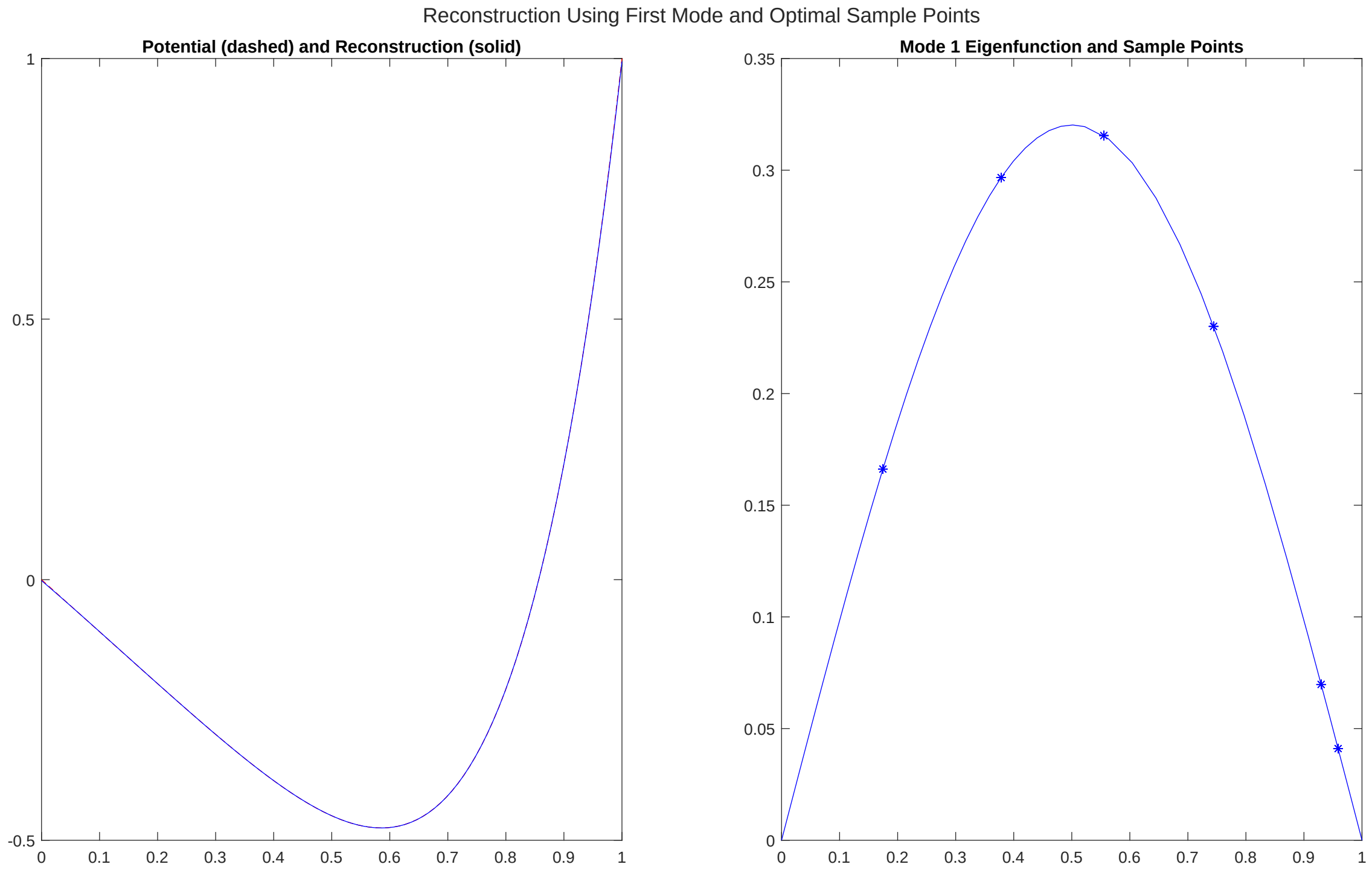

*Example* 7. In the next example, we use the potential $q(x) = 1 + \frac{1}{2}\sin(4\pi x)$. We show a reconstruction using the Legendre basis and one using the basis $\{\cos 2(l-1)\pi x, \sin 2l\pi x\}_{l=1}^{2}$. Clearly the second basis should be the better (and should recover $q$ exactly). But suppose we don't know this.

The reconstruction with the Legendre basis is:

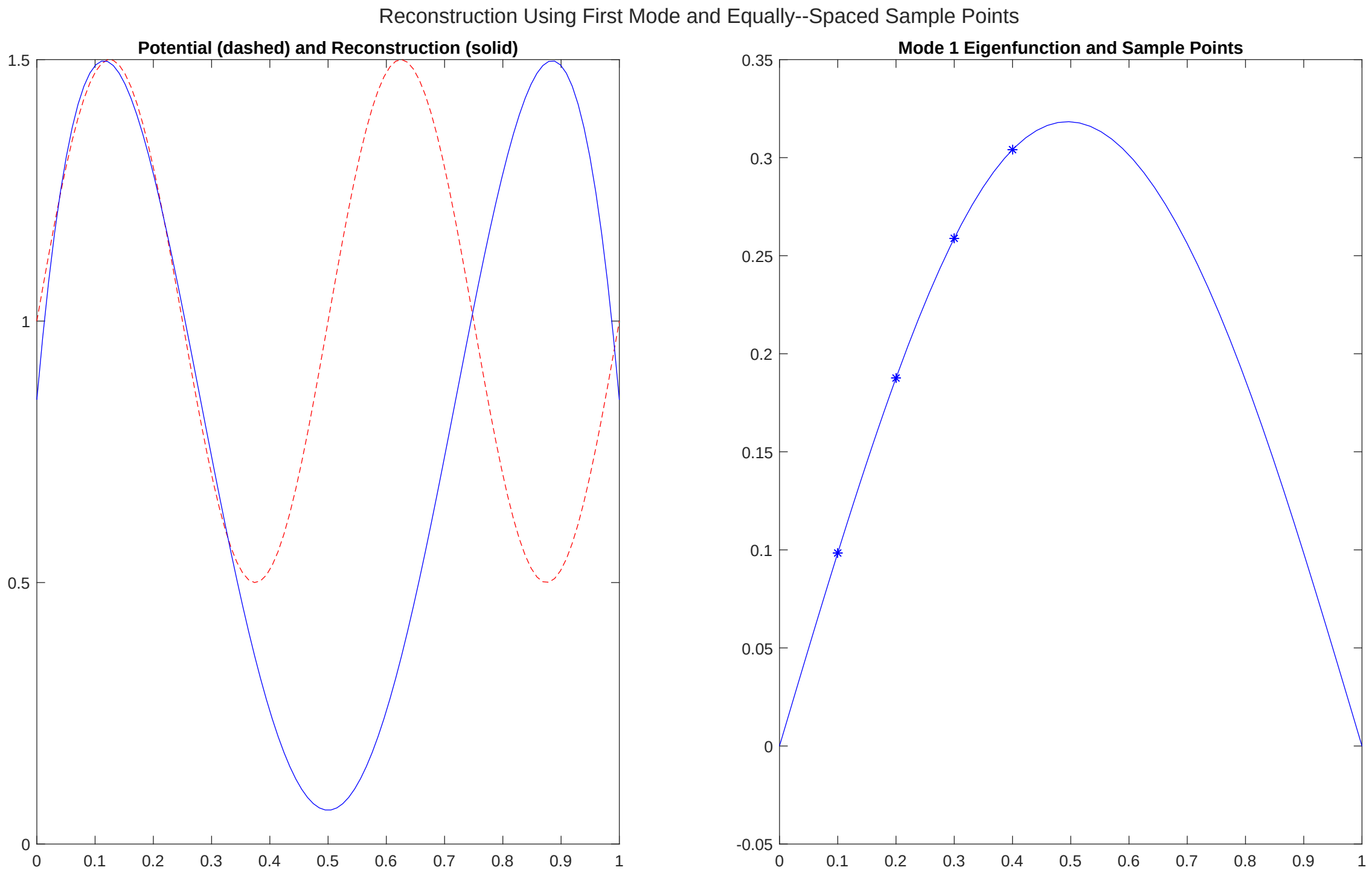

The reconstruction with the trigonometric basis is:

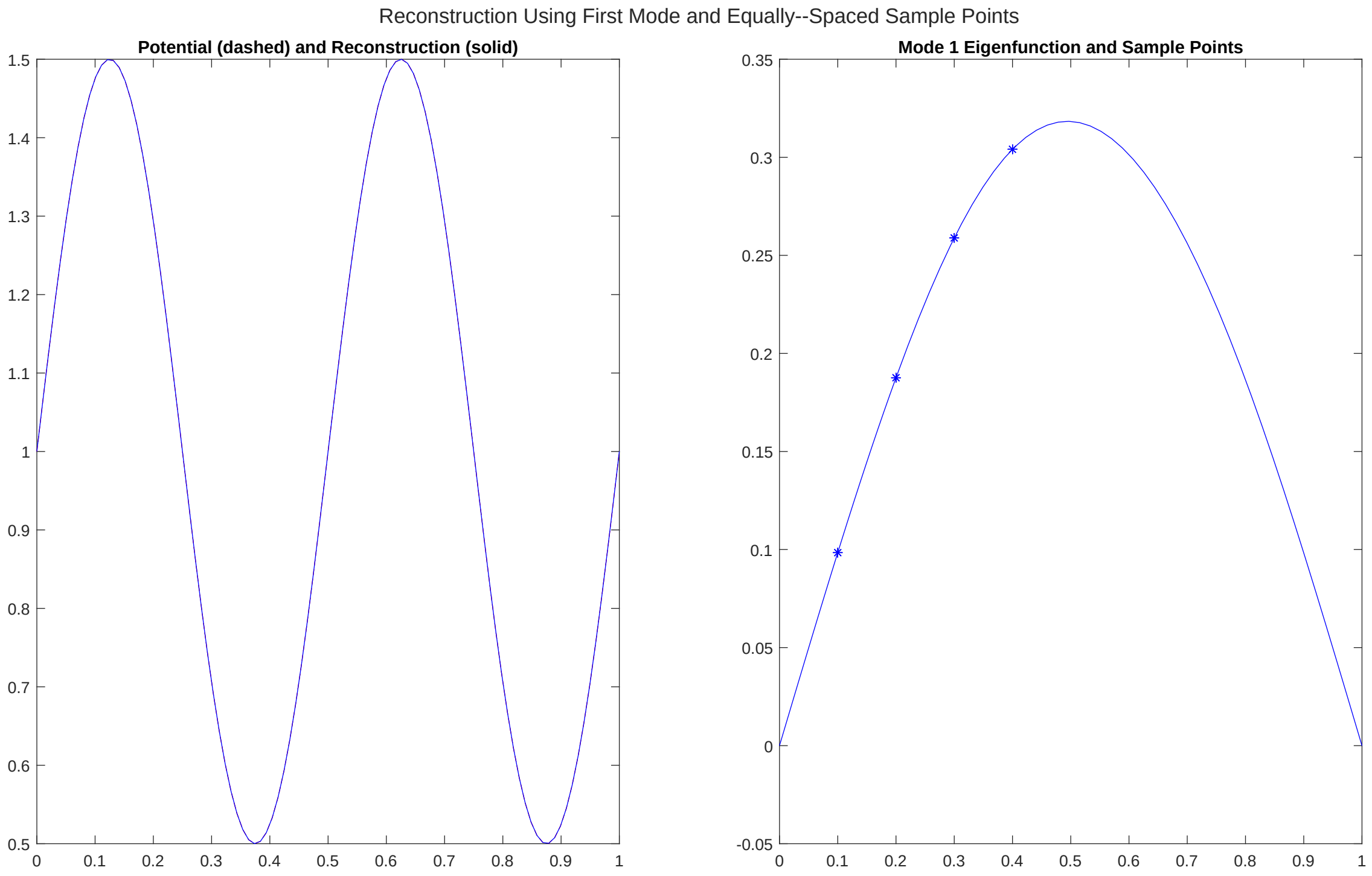

The $L^2$ norm of the potential recovered using the Legendre basis is $1.7306$ and the one using the trigonometric basis is $2$. Therefore, our selection criteria indicates that the trigonometric basis gives the better approximation (which is true in this case, of course).

*Example* 8. In the final exam, we begin with a Legendre basis and a trigonometric basis. Then we pick the four best basis functions (based on the criteria selection mentioned above) and get a better approximation (in this case, it is exact). The target potential is $q(x) = 1 + (t - .5)^2 + \frac{1}{2}\sin(4\pi x)$. In the Legendre basis, the approximation is:

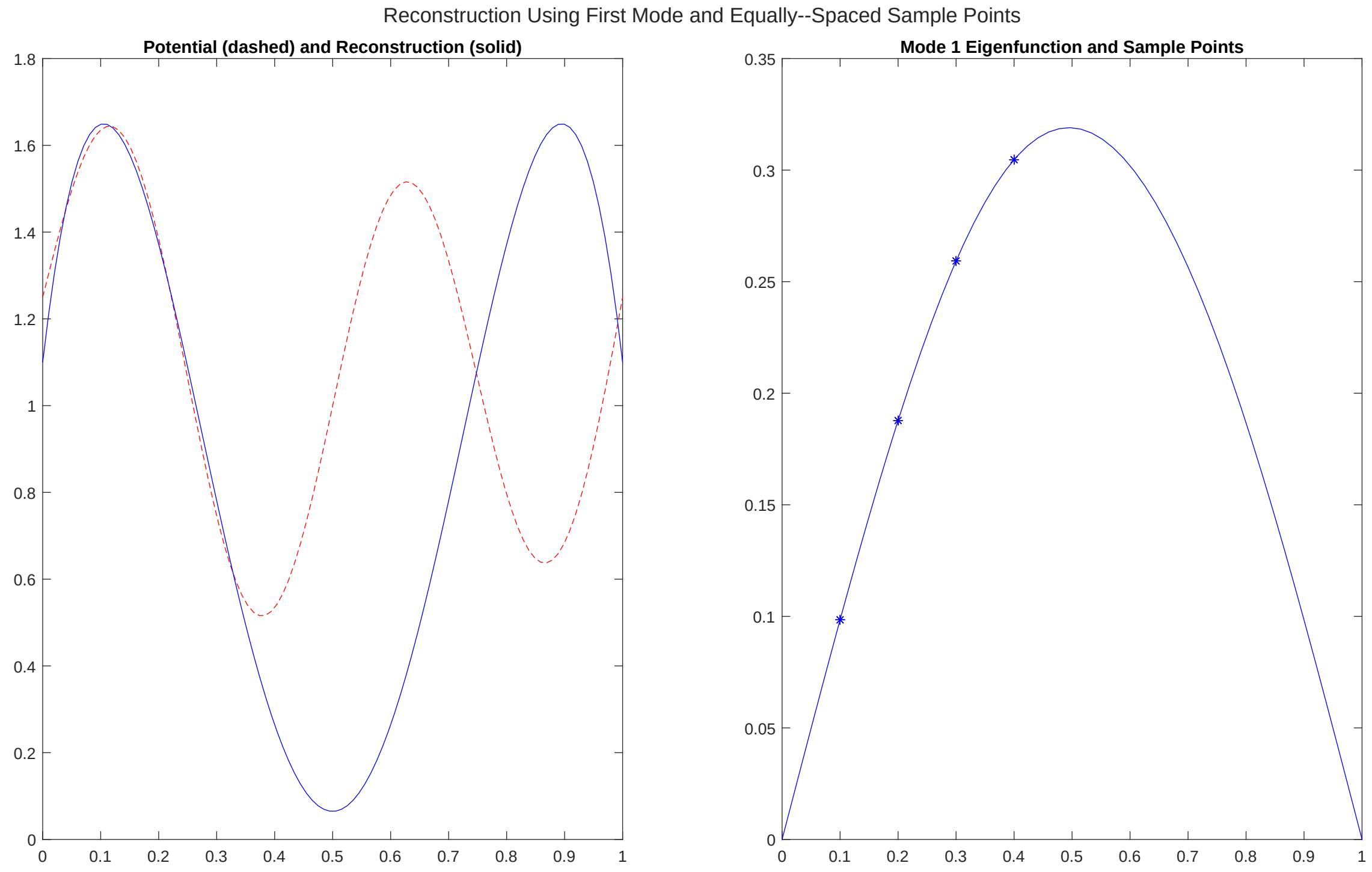

And in the trigonometric basis, the approximation is:

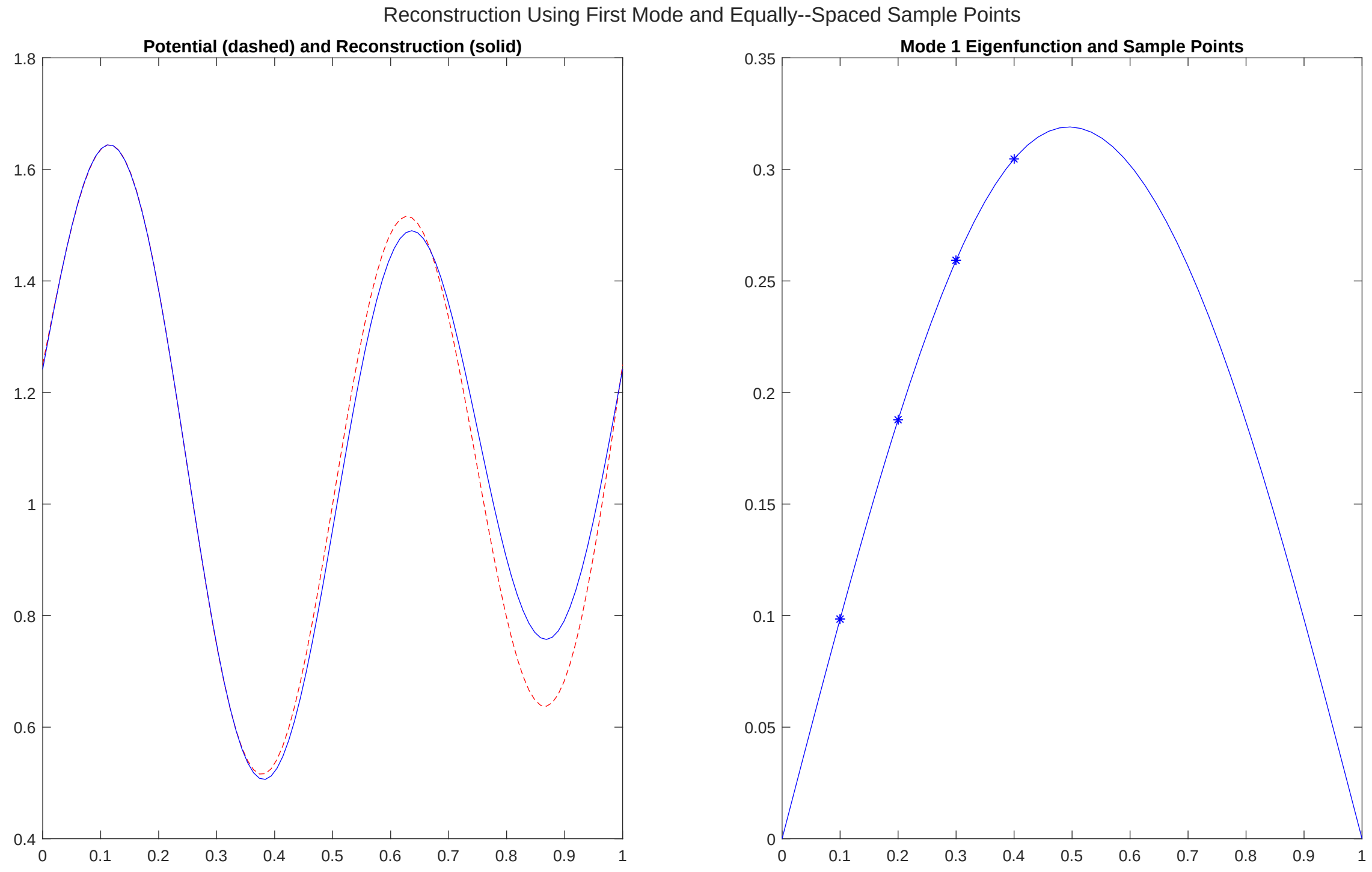

By looking at the coefficients, we can pick the four largest which correspond to the constant function from both bases, the quadratic Legendre polynomial, and the $\sin(4\pi x)$ function from the trigonometric basis. Using these four basis functions the reconstruction is:

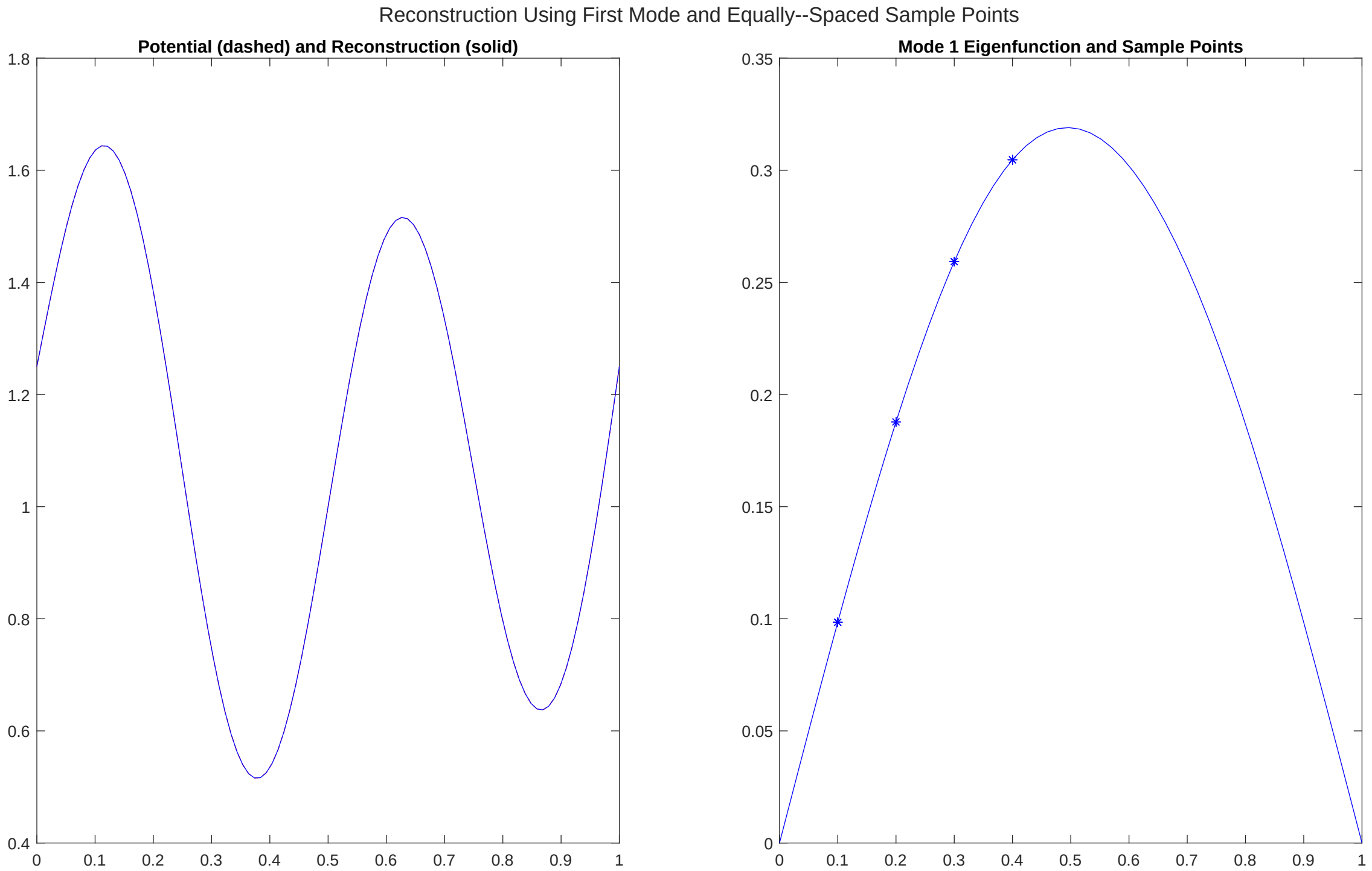

## 4. Appendix

In this section is the code for the inverse problem. Below, the function "jac" computes the matrix in the iteration. As above, the function ivp is just an initial value problem solve.

```
%
% Given N Dirichlet eigenvalues for the problem -y'' + qy = l y
% X is the sample points and L is the data. We need to have
% a Y vector as well; this is the LHS.
% B is the basis functions.
%
% X must be a row vector and contains *only* the sample points
% L must be a row vector and contains the evals
% Y must be a column vector and contains the rhs of the solution

function [coefs, iters] = inverse_full(X, L, Y, B, mode, num_coefs)
    %% estimate qbar
    qbar_est = L(mode) - (mode.*pi).^2;
    L = L - qbar_est;

    mat = jac(X,L,B);
    [U, S, V] = svd(mat);
    d = diag(S);
    CUT_OFF = 1e-6;
    d(abs(d)<CUT_OFF) = 0;
```

```
    S = diag(d);
    mat = U*S*V';
    mati = pinv(mat);

    tolerance = 1e-4;
    iters = 0;
    delta = ones(num_coefs, 1);
    coefs = zeros(num_coefs, 1);
    while ((iters < 10)) && (max(abs(delta)) > tolerance)
        iters = iters + 1;
        rhs = get_rhs(coefs, X, L, Y, B, num_coefs);
        delta = mati*rhs;
        coefs = coefs + delta;
    end

    coefs(1) = coefs(1) + qbar_est;
end


%
% Computes the RHS for the iteration
%
function rhs = get_rhs(coefs, X, L, Y, B, num_coefs)
    tol = [1e-10 1e-10 1e-11];
    N = num_coefs;
    rhs = zeros(N, 1);
    for row = 1:N
       q = @(t) pot(t, coefs, B);
       [~, u] = ivp(q, L(row), [0 X(row)], tol);
       rhs(row) = u(end,1) - Y(row);
    end
end
```

Below is the listing for the "jac" function:

```
%{
This will give the matrix in the iteration for the
Newton's method in the inverse problem.

X is the sample points
L is the eigenvalues
B is the basis
%}

function jac = jac(X, L, B)
```

```
    N = size(X, 2);
    jac = zeros(N, N);

    for row = 1:N
        for col = 1:N
            xpt = X(row);
            jac(row, col) = integral(@(t)jf(t, xpt, B{col}, L(row)), 0, xpt);
        end
    end
end


%
% This will give the function for the entry in the jacobian that is to
% be integrated.
%
function val = jf(t, xpt, b, l)
    val = sin(sqrt(l).*t) .* sin(sqrt(l).*(t-xpt)) .* b(t);
    val = val ./ l;

end
```

College Station, TX
*Email address*: `robrahm@tamu.edu`